\documentclass[12pt,reqno]{amsart}
\usepackage{amsmath,amssymb,amsthm,mathtools,microtype,enumitem,url}
\usepackage[margin=.9in]{geometry}
\usepackage[hidelinks]{hyperref}
\newtheorem{theorem}{Theorem}[section]
\newtheorem{corollary}[theorem]{Corollary}
\newtheorem{proposition}[theorem]{Proposition}
\newtheorem{lemma}[theorem]{Lemma}

\theoremstyle{definition}
\newtheorem{definition}[theorem]{Definition}

\theoremstyle{remark}
\newtheorem{remark}[theorem]{Remark}
\DeclareMathOperator{\supp}{supp}

\newcommand{\Cl}{\operatorname{Cl}}
\newcommand{\cO}{\mathcal O}

\title[Support-sensitive zero-sum bounds]{Support-sensitive bounds for shortest zero-sum subsequences}

\author{Claudiu G. Pop}
\address{Mathematical Institute, Leiden University, Einsteinweg 55, 2333 CC Leiden, The Netherlands}
\email{popclaudiu38@gmail.com}

\author{George C. \c{T}urca\c{s}}
\address{Babe\c{s}-Bolyai University, Cluj-Napoca, Romania}
\email{george.turcas@ubbcluj.ro}

\hypersetup{
  pdftitle={Support-sensitive bounds for shortest zero-sum subsequences},
  pdfauthor={Claudiu G. Pop and George C. Turcas}
}

\begin{document}

\begin{abstract}
For a sequence $S$ over a finite abelian group, let $MZ(S)$ denote the length of the shortest nonempty zero-sum subsequence of $S$. We prove that if $G$ is finite abelian of order $n$ and $S$ has length $n$, then $MZ(S)\le n-|\supp(S)|+1$. The same bound holds for every sequence of length at least $|G|$. In cyclic groups we combine this elementary support bound with the Savchev--Chen structure theorem for long zero-sumfree sequences and obtain the sharper estimate $MZ(S)\le n-t(t-1)/2$, where $t=|\supp(S)|$, whenever $S$ has length $n$ over $C_n$ and $MZ(S)-1>n/2$. As a consequence, every length-$n$ sequence over $C_n$ with support size $3$ has a zero-sum subsequence of length at most $n-3$, and this is sharp for $n\ge 5$. We also give an arithmetic application to products of prime ideals in a number field, phrased in the standard class-group and block-monoid setting and a corresponding cyclic class-group sharpening.
\end{abstract}

\dedicatory{Dedicated to Professor Dorin Andrica, on his 70th birthday, in gratitude for his encouragement of young mathematicians and for the many inspiring questions he shared with generations of students.}

\maketitle

\section{Introduction}

Zero-sum theory often asks how long a sequence over a finite abelian group must be before it is forced to contain a zero-sum subsequence.  The Davenport constant is the classical example, and it remains central both in additive combinatorics and in factorization theory; see, for example, the survey of Gao and Geroldinger \cite{GaoGer06} and the monograph of Geroldinger and Halter-Koch \cite{GeroldingerHalterKoch2006}.  In this paper we look at a slightly more local question.  Once a sequence is already in front of us, how does the number of distinct group elements appearing in it constrain the length of its shortest zero-sum subsequence?

To make this precise, for a sequence $S$ over a finite abelian group we write $MZ(S)$ for the minimum length of a nonempty zero-sum subsequence of $S$.  The support $\supp(S)$ is the set of elements occurring in $S$.  Our first result is the support-sensitive bound
\begin{equation} \label{eq:thm31}
  MZ(S)\le |G|-|\supp(S)|+1
\end{equation}
for every length-$|G|$ sequence over a finite abelian group $G$.  The proof is short: remove one occurrence of many distinct support elements, apply the elementary lower bound $|\Sigma(P)|\ge |P|$ to the remaining zero-sumfree sequence, and count inside $G\setminus\{0\}$.  The same argument also gives the natural length-at-least-$|G|$ corollary by passing to a support-preserving subsequence of length $|G|$.

The point of view is close to the support-sensitive literature on subsequence sums.  Eggleton and Erd\H{o}s \cite{EggletonErdos1972} gave early lower bounds for the set of subsequence sums of zero-sumfree sequences, and Freeze and Smith \cite{Freeze2001} proved the sharper estimate
\[
  |\Sigma(P)|\ge |P|+|\supp(P)|-1
\]
for zero-sumfree sequences.  More recent work, such as Wang--Chao--Peng \cite{WangChaoPeng2024}, develops this support-sensitive direction further.  Our main bound uses only the elementary part of this circle of ideas, but it is naturally read as a sequence-level companion to those estimates: larger support leaves less room for all zero-sums to be long.

The second part of the paper keeps the group cyclic and uses genuinely cyclic structure.  Savchev and Chen proved that every sufficiently long zero-sumfree sequence in $C_n$ can be put, after multiplication by a unit, into a form whose least positive representatives have total sum less than $n$ \cite{SavchevChen2006}.  Combining this with the finite-abelian support bound gives
\begin{equation} \label{eq:thm41}
  MZ(S)\le n-\binom{|\supp(S)|}{2}
\end{equation}
whenever $S$ has length $n$ over $C_n$ and $MZ(S)-1>n/2$.  In particular, the support-three case has the exact extremal value $n-3$ for $n\ge 5$, witnessed by the sequence $g^{[n-2]}(2g)(3g)$.

Finally we record an arithmetic application.  We formulate it for products of prime ideals, which is the setting needed for the direct ideal-factorization argument, or equivalently through the usual transfer from a Krull monoid to the block monoid over the class group.  If $K$ is a number field with class group of order $h$, and if $m\ge h$ prime ideals occupy $t$ ideal classes, then some subproduct of at most $h-t+1$ of them is principal and generated by an irreducible element.  When the class group is cyclic, the cyclic zero-sum refinements above give a stronger bound in the same prime-ideal setting.  The restriction to prime ideals is essential for the direct ideal-factorization proof.

\section{Preliminaries}

Throughout, $G$ is a finite abelian group written additively, and $C_n$ denotes the cyclic group of order $n$.  A sequence $S$ over $G$ is a finite multiset of elements of $G$, and subsequences are taken as submultisets of occurrences.  We write $T\mid S$ when $T$ is a subsequence of $S$, $|S|$ for the length of $S$, $\supp(S)$ for the set of elements appearing in $S$, and $h(S)$ for the maximum multiplicity of an element in $S$.  Thus $|\supp(S)|$ is the support size.  A sequence is zero-sumfree if it has no nonempty zero-sum subsequence.  For $g\in G$ and $k\ge 0$, the notation $g^{[k]}$ denotes the sequence consisting of $k$ copies of $g$, with $g^{[0]}$ the empty sequence.

\begin{definition}\label{def:subsequence-sums}
Let $S$ be a sequence over $G$, and let $\sigma(T)$ denote the sum of the terms of a subsequence $T\mid S$.  We set
\[
  \Sigma(S)=\{\sigma(T):T\mid S,\ |T|\ge 1\},
\]
the set of sums of nonempty subsequences of $S$.  For an integer $\ell\ge 0$, let
\[
  \Sigma_\ell(S)=\{\sigma(T):T\mid S,\ |T|=\ell\}.
\]
For $\ell\ge 1$, define
\[
  \Sigma_{\leq \ell}(S)=\bigcup_{1\le j\le \ell}\Sigma_j(S),
  \qquad
  \Sigma_{\geq \ell}(S)=\bigcup_{j\ge \ell}\Sigma_j(S).
\]
In particular, $\Sigma(S)=\Sigma_{\geq 1}(S)$.
\end{definition}

\begin{definition}\label{def:shortest-zero-sum-length}
For a sequence $S$ over $G$, define $MZ(S)$ to be the smallest positive integer $r$ such that $S$ has a zero-sum subsequence of length $r$.  If no nonempty zero-sum subsequence exists, put $MZ(S)=\infty$.
\end{definition}

\begin{lemma}[Zero-sumfree sumset growth]\label{lem:sumset-growth}
Let $P$ be a zero-sumfree sequence over an abelian group.  Then
\[
  |\Sigma(P)|\ge |P|.
\]
\end{lemma}

\begin{proof}
If $P$ is empty, the assertion is trivial.  Otherwise write
\[
  P=g_1\cdots g_k
\]
in an arbitrary order.  For $1\le j\le k$, put
\[
  s_j=g_1+\cdots+g_j.
\]
Each $s_j$ belongs to $\Sigma(P)$, since it is the sum of the nonempty subsequence $g_1\cdots g_j$.

We claim that $s_1,\ldots,s_k$ are pairwise distinct.  If $s_i=s_j$ for some $i<j$, then
\[
  g_{i+1}+\cdots+g_j=0,
\]
which is the sum of a nonempty subsequence of $P$, contradicting zero-sumfreeness.  Hence $\Sigma(P)$ contains at least the $k$ distinct elements $s_1,\ldots,s_k$, and so the conclusion follows.
\end{proof}

\begin{remark}\label{rem:freeze-smith-bound}
Freeze--Smith prove the stronger support-sensitive estimate
\[
  |\Sigma(P)|\ge |P|+|\supp(P)|-1
\]
for zero-sumfree sequences.  We only use the elementary bound above in the proof of the main theorem.
\end{remark}

\section{The finite-abelian support bound}

We begin with the group-theoretic bound.  No cyclic ordering, generator, or residue choice is used here; only the size of the ambient group enters the proof.

\begin{theorem}\label{thm:finite-abelian-main}
Let $G$ be a finite abelian group of order $n$, and let $S$ be a sequence over $G$ with $|S|=n$.  Then
\[
  MZ(S)\le n-|\supp(S)|+1.
\]
Equivalently, if $MZ(S)=n-s$ for some $s\in\{0,\ldots,n-1\}$, then $|\supp(S)|\le s+1$.
\end{theorem}

\begin{proof}
Put $t=|\supp(S)|$.  First, $S$ cannot be zero-sumfree.  Indeed, if $S$ were zero-sumfree, Lemma~\ref{lem:sumset-growth} would give $|\Sigma(S)|\ge |S|=n$, while zero-sumfreeness and the convention that $\Sigma(S)$ contains only nonempty subsequence sums would give $\Sigma(S)\subseteq G\setminus\{0\}$, whose cardinality is $n-1$.  This is impossible.  Thus $MZ(S)$ is finite; write $r=MZ(S)$, so $1\le r\le n$.

If $r=1$, then the desired inequality is $1\le n-t+1$, which follows from $t\le n$.  If $r=2$, then $0$ is not a term of $S$, for otherwise $r=1$.  Hence $\supp(S)\subseteq G\setminus\{0\}$, so $t\le n-1$, and again $r\le n-t+1$.

Assume now that $r\ge 3$ and suppose for contradiction that $t\ge n-r+2$.  Choose distinct support elements $a_1,\ldots,a_{n-r+2}$ occurring in $S$, and remove one occurrence of each.  The remaining sequence $P$ has length $r-2$.  Since $|P|<r$, it is zero-sumfree.  By Lemma~\ref{lem:sumset-growth}, $|\Sigma(P)|\ge r-2$.

Let $A=\{-a_1,\ldots,-a_{n-r+2}\}$.  Since $r>1$, zero is not a term of $S$, and hence $A\subseteq G\setminus\{0\}$.  The map $x\mapsto -x$ is a bijection on $G$, so $|A|=n-r+2$.  Also $\Sigma(P)\subseteq G\setminus\{0\}$.  But
\[
 |\Sigma(P)|+|A|\ge (r-2)+(n-r+2)=n>|G\setminus\{0\}|.
\]
Therefore $\Sigma(P)\cap A\ne \varnothing$.  Choose $x\in\Sigma(P)\cap A$, say $x=-a_i$.  By the definition of $\Sigma(P)$, there is a nonempty subsequence $Q\mid P$ with $\sigma(Q)=x$.  Adding the selected occurrence of $a_i$, which is not an occurrence of $P$, gives a zero-sum subsequence of $S$ of length at most
\[
  |P|+1=r-1,
\]
contradicting the definition of $r$.  Thus $t\le n-r+1$, which is equivalent to $r\le n-t+1$.
\end{proof}

We now present a corollary for the special case in which the length of the sequence is at least $|G|$.

\begin{corollary}\label{cor:length-at-least}
Let $G$ be a finite abelian group of order $n$.  If $S$ is a sequence over $G$ with $|S|\ge n$, then
\[
  MZ(S)\le n-|\supp(S)|+1.
\]
\end{corollary}

\begin{proof}
Put $t=|\supp(S)|$.  Since $\supp(S)\subseteq G$, we have $t\le n$.  Choose one occurrence in $S$ of each element of $\supp(S)$, and call the resulting subsequence $U$.  Then $|U|=t$ and $\supp(U)=\supp(S)$.  After removing $U$ from $S$, at least $n-t$ terms remain, since $|S|\ge n$.  Choose any $n-t$ of those remaining terms, appending no terms if $n-t=0$, and add them to $U$.  This gives an $n$-term subsequence $T\mid S$ with $\supp(T)=\supp(S)$.

By Theorem~\ref{thm:finite-abelian-main}, $T$ has a nonempty zero-sum subsequence of length at most
\[
  n-|\supp(T)|+1=n-|\supp(S)|+1.
\]
The same subsequence is a subsequence of $S$, proving the claim.
\end{proof}

\section{Cyclic refinements via Savchev--Chen}

The finite-abelian theorem is linear in the support size.  In cyclic groups, long zero-sumfree sequences have much more structure, and this lets us improve the bound in the range where the shortest zero-sum subsequence is itself long.  The finite-abelian theorem is used below only to obtain the preliminary numerical bound on the support; the refinement is cyclic because it relies on Savchev--Chen's structure theorem and on least positive representatives in $C_n$.

\begin{theorem}[Cyclic long-regime refinement]\label{thm:savchev-chen-refinement}
Let $S$ be a length-$n$ sequence over $C_n$.  Put $r=MZ(S)$ and $t=|\supp(S)|$.  If $r-1>n/2$, then
\[
  r\le n-\binom{t}{2}.
\]
\end{theorem}

\begin{proof}
The cases $n\le 2$ are vacuous, since Theorem~\ref{thm:finite-abelian-main} gives $r\le n$ and hence $r-1\le n/2$.  Assume $n\ge 3$.  Again by Theorem~\ref{thm:finite-abelian-main},
\[
  t\le n-r+1.
\]
Since $r-1>n/2$, we have $n-r+1<r-1$, and hence $t\le r-1$.

Choose one occurrence of each support element of $S$ and extend this choice arbitrarily to a subsequence $T\mid S$ of length $r-1$.  This is possible because $t\le r-1\le n$.  By construction, $\supp(T)=\supp(S)$, so $T$ has exactly $t$ distinct terms.  Moreover, $T$ is zero-sumfree: otherwise a nonempty zero-sum subsequence of $T$ would be a nonempty zero-sum subsequence of $S$ of length at most $r-1$, contradicting the definition of $r=MZ(S)$.

Since $|T|=r-1>n/2$, Savchev--Chen's theorem \cite[Theorem~8]{SavchevChen2006} applies.  Thus there is a unit $u$ modulo $n$ such that the sum of the least positive representatives of the terms of $uT$ is less than $n$.  Multiplication by $u$ is an automorphism of $C_n$, so it preserves zero-sumfreeness and the number of distinct terms.  Since $T$ is zero-sumfree, neither $T$ nor $uT$ contains $0$; hence all these least positive representatives lie in $\{1,\ldots,n-1\}$.

It remains only to estimate a positive integer multiset.  A multiset of $r-1$ positive integers with exactly $t$ distinct values has total sum at least
\[
  (r-1)+\binom{t}{2}.
\]
Indeed, if the distinct values are $1\le v_1<\cdots<v_t$, then $v_i\ge i$ for every $i$, and the remaining $r-1-t$ terms are each at least $1$.  Thus the total sum is at least
\[
  (1+\cdots+t)+(r-1-t)=(r-1)+\binom{t}{2}.
\]
The Savchev--Chen representative sum is an integer strictly less than $n$, so it is at most $n-1$.  Therefore
\[
  (r-1)+\binom{t}{2}\le n-1,
\]
which rearranges to $r\le n-\binom{t}{2}$.
\end{proof}

In the special case in which $|\supp(S)| =3$, this yields the following corollary.

\begin{corollary}\label{cor:support-three-cyclic}
Let $n\ge 5$, and let $S$ be a length-$n$ sequence over $C_n$ with $|\supp(S)|=3$.  Then
\[
  MZ(S)\le n-3.
\]
This bound is sharp: if $g$ is a generator of $C_n$, then
\[
  MZ\bigl(g^{[n-2]}(2g)(3g)\bigr)=n-3.
\]
\end{corollary}

\begin{proof}
Put $r=MZ(S)$.  By Theorem~\ref{thm:finite-abelian-main},
\[
  r\le n-|\supp(S)|+1=n-2.
\]
Assume first that $n\ge 7$.  If $r>n-3$, then the preceding inequality forces $r=n-2$.  Hence
\[
  r-1=n-3>n/2,
\]
so Theorem~\ref{thm:savchev-chen-refinement}, applied with $|\supp(S)|=3$, gives
\[
  r\le n-\binom{3}{2}=n-3,
\]
a contradiction.  Thus $r\le n-3$ for $n\ge 7$.

It remains to handle $n=5$ and $n=6$.  For $n=5$, if $0\in\supp(S)$, then $MZ(S)=1\le 2=n-3$.  Otherwise $\supp(S)$ is a three-element subset of the four nonzero elements of $C_5$.  These four elements split into two inverse pairs, so any three-element subset contains one inverse pair.  The corresponding two occurrences form a zero-sum subsequence, and $MZ(S)\le 2=n-3$.

For $n=6$, write $C_6=\langle g\rangle$ and identify terms with their coefficients modulo $6$.  If $0\in\supp(S)$, then $MZ(S)=1\le 3=n-3$.  If $\supp(S)$ contains one of the inverse pairs $\{g,5g\}$ or $\{2g,4g\}$, then $S$ has a zero-sum subsequence of length $2$.

We may therefore assume that $\supp(S)$ contains no zero and no inverse pair.  Then it must contain $3g$ and one element from each of the pairs $\{g,5g\}$ and $\{2g,4g\}$.  Thus the support, written in coefficients of $g$, is one of
\[
  \{1,2,3\},\quad \{1,3,4\},\quad \{2,3,5\},\quad \{3,4,5\}.
\]
Let $m_a$ denote the multiplicity of $ag$ in $S$.  If $m_3\ge 2$, then $3g+3g=0$, so $MZ(S)\le 2$.  We may assume $m_3=1$.

For support $\{1,2,3\}$, the subsequence $g(2g)(3g)$ has sum $0$.  For support $\{3,4,5\}$, the subsequence $(3g)(4g)(5g)$ has sum $0$.  For support $\{1,3,4\}$, we have $m_1+m_4=5$.  If $m_1\ge 2$, then $g+g+4g=0$; if $m_1\le 1$, then $m_4\ge 4$, and $4g+4g+4g=0$.  For support $\{2,3,5\}$, we have $m_2+m_5=5$.  If $m_2\ge 3$, then $2g+2g+2g=0$; if $m_2\le 2$, then $m_5\ge 3$, and $2g+5g+5g=0$.  Hence $MZ(S)\le 3=n-3$ in every case.

It remains to prove sharpness.  Let
\[
  S_n=g^{[n-2]}(2g)(3g)
\]
with $g$ a generator of $C_n$.  Since $n\ge 5$, the elements $g,2g,3g$ are distinct.  The subsequence $g^{[n-5]}(2g)(3g)$ has length $n-3$ and sum $ng=0$, so $MZ(S_n)\le n-3$.

Now consider any nonempty subsequence of $S_n$.  It has the form
\[
  g^{[a]}(2g)^{[\epsilon_2]}(3g)^{[\epsilon_3]},
\]
where $0\le a\le n-2$ and $\epsilon_2,\epsilon_3\in\{0,1\}$.  Its length is
\[
  L=a+\epsilon_2+\epsilon_3,
\]
and its coefficient of $g$ is
\[
  c=a+2\epsilon_2+3\epsilon_3=L+\epsilon_2+2\epsilon_3.
\]
If $L<n-3$, then $L\le n-4$.  Since the subsequence is nonempty,
\[
  1\le c\le L+3\le n-1.
\]
Thus $c$ is not divisible by $n$, and the subsequence is not zero-sum.  Therefore no zero-sum subsequence has length less than $n-3$, while one of length $n-3$ exists.  Hence $MZ(S_n)=n-3$.
\end{proof}

\section{Arithmetic applications}

Given a number field $K$ with ring of integers $\mathcal O_K$, as mentioned previously in the introduction, it is known that the Davenport constant of the class group $Cl(K)$ of $\mathcal O_K$ is related to its \textit{elasticity} (see \cite{Pollack2025}), i.e. to lengths of distinct irreducible factorizations of nonzero elements in $\mathcal O_K$. Motivated by this connection, we now record the arithmetic form of Corollary~\ref{cor:length-at-least}.  The statement is deliberately made for prime ideals.  For arbitrary integral ideals the corresponding subproduct assertion need not hold: a principal factor of a product of arbitrary integral ideals need not be a product of a subcollection of the given ideals.

\begin{proposition}\label{prop:prime-ideal-application}
Let $K$ be a number field with ring of integers $\cO_K$, and let $G=\Cl(K)$ have order $h$.  Let $\mathfrak p_1,\ldots,\mathfrak p_m$ be nonzero prime ideals of $\cO_K$ with $m\ge h$, and suppose their ideal classes occupy $t$ distinct classes of $G$.  Then there is a nonempty index set $I\subseteq\{1,\ldots,m\}$ such that
\[
 |I|\le h-t+1
\]
and
\[
 \prod_{i\in I}\mathfrak p_i=(\pi)
\]
for some irreducible element $\pi\in\cO_K$.  Repeated prime ideals are counted as distinct indexed occurrences.
\end{proposition}

\begin{proof}
Form the indexed sequence of ideal classes
\[
  S=[\mathfrak p_1]\cdots[\mathfrak p_m]
\]
over $G=\Cl(K)$, written additively.  This sequence has length $m\ge h$ and support size $t$.  By Corollary~\ref{cor:length-at-least},
\[
  MZ(S)\le h-t+1.
\]
Choose a nonempty index set $I\subseteq\{1,\ldots,m\}$ with $|I|=MZ(S)$ and
\[
  \sum_{i\in I}[\mathfrak p_i]=0.
\]
Then
\[
  \prod_{i\in I}\mathfrak p_i=(\pi)
\]
for some nonzero $\pi\in\cO_K$.  Since $I$ is nonempty and each $\mathfrak p_i$ is a proper ideal, this product is a proper ideal, so $\pi$ is not a unit.

We prove that $\pi$ is irreducible.  Suppose, to the contrary, that $\pi=ab$ with nonunits $a,b\in\cO_K$.  Then
\[
  (a)(b)=(\pi)=\prod_{i\in I}\mathfrak p_i.
\]
For each nonzero prime ideal $\mathfrak q$, put
\[
  e_{\mathfrak q}=\#\{i\in I:\mathfrak p_i=\mathfrak q\}.
\]
Unique factorization of nonzero ideals into prime ideals gives
\[
  v_{\mathfrak q}((a))+v_{\mathfrak q}((b))=e_{\mathfrak q}
\]
for every $\mathfrak q$, and no prime outside the indexed multiset $\{\mathfrak p_i:i\in I\}$ occurs in $(a)$ or $(b)$.  Since $a$ is a nonunit, $(a)$ is a proper nonzero integral ideal, so its prime-ideal factorization is nonempty.  Since $b$ is a nonunit, $(b)$ is also proper, so $(a)$ cannot use all indexed prime-ideal occurrences.  Indeed, if $v_{\mathfrak q}((a))=e_{\mathfrak q}$ for every $\mathfrak q$, then all valuations of $(b)$ would be zero, so $(b)=\cO_K$, contradicting that $b$ is a nonunit.

Choose, for each $\mathfrak q$, exactly $v_{\mathfrak q}((a))$ of the indices $i\in I$ with $\mathfrak p_i=\mathfrak q$.  This gives a nonempty proper subset $J\subsetneq I$ such that
\[
  (a)=\prod_{j\in J}\mathfrak p_j.
\]
Since $(a)$ is principal, we get
\[
  \sum_{j\in J}[\mathfrak p_j]=0
\]
in $\Cl(K)$.  This is a nonempty zero-sum subsequence of $S$ of length $|J|<|I|=MZ(S)$, a contradiction.  Hence no such factorization $\pi=ab$ exists, and $\pi$ is irreducible.
\end{proof}

When the class group is cyclic, Proposition~\ref{prop:prime-ideal-application} inherits the cyclic zero-sum refinements from the preceding section.  The resulting bound is best read as a dichotomy: either the shortest principal subproduct is already no longer than about half the class number, or the long zero-sumfree structure theorem applies.

\begin{proposition}[Cyclic class-group sharpening]\label{prop:cyclic-class-group-sharpening}
Let $K$ be a number field with ring of integers $\cO_K$, and suppose that $\Cl(K)$ is cyclic of order $h$.  Let $\mathfrak p_1,\ldots,\mathfrak p_m$ be nonzero prime ideals of $\cO_K$ with $m\ge h$, and suppose their ideal classes occupy $t$ distinct classes.  Put
\[
 B_h(t)=
 \min\left\{
 h-t+1,\,
 \max\left\{\left\lfloor\frac h2\right\rfloor+1,\,
 h-\binom{t}{2}\right\}
 \right\}.
\]
Then there is a nonempty index set $I\subseteq\{1,\ldots,m\}$ such that
\[
 |I|\le B_h(t)
\]
and
\[
 \prod_{i\in I}\mathfrak p_i=(\pi)
\]
for some irreducible element $\pi\in\cO_K$.  If, in addition, $h\ge 5$ and $t=3$, then one can take $|I|\le h-3$.
\end{proposition}

\begin{proof}
Form the indexed class sequence
\[
  S=[\mathfrak p_1]\cdots[\mathfrak p_m]
\]
over $\Cl(K)$.  Choose one occurrence from each of the $t$ represented ideal classes and then add arbitrary further occurrences until the resulting subsequence $T\mid S$ has length $h$.  This is possible because $m\ge h$, and it gives
\[
  |T|=h,\qquad |\supp(T)|=t.
\]
Identify $\Cl(K)$ with $C_h$, and put $q=MZ(T)$.

By Corollary~\ref{cor:length-at-least}, or by Theorem~\ref{thm:finite-abelian-main} applied directly to $T$,
\[
  q\le h-t+1.
\]
On the other hand, either $q-1\le h/2$, in which case
\[
  q\le \left\lfloor\frac h2\right\rfloor+1,
\]
or $q-1>h/2$, in which case Theorem~\ref{thm:savchev-chen-refinement} gives
\[
  q\le h-\binom{t}{2}.
\]
Combining the linear bound with this cyclic dichotomy gives $q\le B_h(t)$.  Since $T$ is a subsequence of the full class sequence $S$, the full sequence also has a nonempty zero-sum subsequence of length at most $B_h(t)$.

Choose a globally shortest nonempty zero-sum subsequence of $S$, and let $I$ be its index set.  Then $|I|\le B_h(t)$ and
\[
  \sum_{i\in I}[\mathfrak p_i]=0,
\]
so
\[
  \prod_{i\in I}\mathfrak p_i=(\pi)
\]
for some nonzero $\pi\in\cO_K$.  The same indexed-prime argument used in Proposition~\ref{prop:prime-ideal-application} shows that $\pi$ is irreducible: a nontrivial factorization $\pi=ab$ would make $(a)$ a principal product of a nonempty proper indexed subcollection of the primes in $I$, giving a shorter nonempty zero-sum subsequence of the class sequence.

Finally assume $h\ge 5$ and $t=3$.  The support-preserving subsequence $T$ has length $h$ and support size $3$ over $C_h$, so Corollary~\ref{cor:support-three-cyclic} gives $MZ(T)\le h-3$.  Repeating the preceding global-minimality and irreducibility argument gives an indexed principal irreducible subproduct of length at most $h-3$.
\end{proof}

In factorization-theoretic language, the selected class sequence in each proposition is an atom of the block monoid over $\Cl(K)$, and the principal generator is an atom of $\cO_K^\bullet$.  This is the ring-of-integers instance of the standard transfer from Krull monoids to monoids of zero-sum sequences over the class group; see \cite{GeroldingerHalterKoch2006,GeroldingerKainrath2021,Schmid2004BlockMonoids}.

\section{Future questions}

For a finite abelian group $G$, define
\[
 f_G(t)=\max\{MZ(S): |S|=|G|,\ |\supp(S)|=t\}.
\]
Theorem~\ref{thm:finite-abelian-main} gives $f_G(t)\le |G|-t+1$, while Corollary~\ref{cor:support-three-cyclic} shows that $f_{C_n}(3)=n-3$ for $n\ge 5$.  It would be natural to determine $f_{C_n}(t)$ for larger fixed $t$, and then to ask which parts of the cyclic picture survive for noncyclic groups. Moreover, it is natural to ask for which values of t equality is attained in the bound above. One can show that the extremal cases arise only when $t\in \{1,2,n-1,n\}$, occurring in symmetric pairs.

There is also a stability question behind the linear bound.  The Freeze--Smith estimate suggests that equality in Theorem~\ref{thm:finite-abelian-main} should force a large multiplicity in the original sequence.  We do not use such a statement here, and we leave a precise reviewed stability theorem for future work.

\bibliographystyle{alpha}
\bibliography{Davenport_references}

\end{document}